\documentclass[11 pt]{amsart}
\usepackage[utf8]{inputenc}
\usepackage[dvipsnames,svgnames,x11names,hyperref]{xcolor}
\definecolor{ffmblue}{HTML}{006092}
\usepackage{comment}
\allowdisplaybreaks
\usepackage{csquotes}
\usepackage[hyphens]{url}
\usepackage[colorlinks=true,linkcolor=Maroon,citecolor=Maroon,urlcolor=ffmblue,hypertexnames=false,linktocpage]{hyperref}
\usepackage{bookmark}
\usepackage{amsmath,thmtools,mathtools}
\mathtoolsset{showonlyrefs=true}

\newcounter{mparcnt}

\usepackage{fancyhdr}
\usepackage{esint}
\usepackage{enumerate}

\usepackage{pictexwd,dcpic}
\usepackage{graphicx}
\usepackage{caption}
\swapnumbers
\declaretheorem[name=Theorem,numberwithin=section]{thm}
\declaretheorem[name=Remark,style=remark,sibling=thm]{rem}

\declaretheorem[name=Definition,style=definition,sibling=thm]{defn}

\declaretheorem[name=Assumption,style=definition,sibling=thm]{assum}

\declaretheorem[style=definition,name=Definition,numbered=no]{definition}

\numberwithin{equation}{section}

\newcommand{\del}{\partial}

\newcommand{\pf}[1]{\begin{proof}#1 \end{proof}}

\begin{document}

\title[Stability of Minkowski inequalities in warped spaces]{Stability of Minkowski-type inequalities in certain warped product spaces}

\author{Prachi Sahjwani}
\address{\flushleft\parbox{\linewidth}{Cardiff University\\ School of Mathematics\\ Senghennydd Road\\ Cardiff CF24 4AG\\Wales\\ {\href{mailto:sahjwanip@cardiff.ac.uk}{sahjwanip@cardiff.ac.uk}} }}

\date{\today}


\begin{abstract}
This paper proves quantitative stability estimates for five Minkowski-type inequalities for hypersurfaces in warped product spaces. In each case we show that if a hypersurface nearly achieves equality, it must be geometrically close, in the Hausdorff sense, to a radial slice. The ambient spaces are warped products $(a,b)\times\mathbb{S}^n$ with metric $dr^2+\lambda^2(r)g_{\mathbb{S}^n}$, as well as the Reissner-Nordstr\"om Anti-de Sitter (RN-AdS) and Anti-de Sitter Schwarzschild (AdS-Schwarzschild) manifolds. The proofs combine two ingredients: a quantitative analysis of locally constrained inverse curvature flows, which yields bounds on the traceless second fundamental form in terms of the deficit in the inequality, and a new rigidity theorem for hypersurfaces in locally conformally flat spaces, which converts such bounds into Hausdorff closeness to a radial slice.
\end{abstract}

\maketitle

\section{Introduction}

Minkowski inequalities play a fundamental role in the study of convex bodies and hypersurface geometry, providing lower bounds for curvature integrals in terms of the surface area or enclosed volume of a hypersurface. The classical Minkowski inequality for a closed convex surface $\Sigma$ in $\mathbb{R}^{3}$ states
\begin{equation}
\label{Minkowski surface}
\sqrt{16\pi\lvert\Sigma\rvert} \leq \int_{\Sigma} H\, d\mu,
\end{equation}
where $H$ is the mean curvature and $\lvert\Sigma\rvert$ is the area of $\Sigma$. In higher dimensions, for a convex hypersurface $\Sigma$ in $\mathbb{R}^{n}$,
\begin{equation}
\label{Minkowski dim n}
(n-1)|\mathbb{S}^{n-1}|^{\frac{1}{n-1}} |\Sigma|^{\frac{n-2}{n-1}} \leq \int_{\Sigma}H\, d\mu.
\end{equation}
Equality in both cases holds if and only if $\Sigma$ is a sphere. It is then natural to ask: how close must $\Sigma$ be to a sphere, given that the deficit in \eqref{Minkowski surface} or \eqref{Minkowski dim n} is small? Stability results for various geometric inequalities have been widely investigated; for the isoperimetric inequality, see \cite{figalli_mass_2010, fusco_sharp_2008, ivaki_stability_2014}, and for quermassintegral inequalities, see \cite{gao_geometric_2024, sahjwani_stability_2023, scheuer_stability_2023, schneider_stability_1990, schneider_convex_2013, zhou_stability_2023}. Extensions of Minkowski-type inequalities to warped product spaces and spaces of non-constant curvature have received significant attention \cite{brendle_minkowski_2016,chen_penrose_2019, natario_minkowski-type_2015, scheuer_minkowski_2021, scheuer_minkowski_2022, wang_minkowski-type_2015, wei_minkowski-type_2018}. The present paper addresses the corresponding stability question in this broader setting.

\subsection*{Main results}

We prove quantitative stability estimates for five Minkowski-type inequalities. The ambient spaces are warped products $(a,b)\times\mathbb{S}^n$, as well as the RN-AdS and AdS-Schwarzschild manifolds. A key feature of all these spaces is that they are locally conformally flat, which plays a central role in our approach.

\begin{rem}\label{remark conformal}
The warped product spaces $M^{n+1} = (a,b) \times \mathbb{S}^n$ with $\bar{g}=dr^2 + \lambda^2(r) g_{\mathbb{S}^n}$ are locally conformally flat: in local coordinates, $\bar{g}_{\alpha\beta} = e^{2\omega} \tilde{g}_{\alpha\beta}$ where $\tilde{g}$ is the Euclidean metric and $\omega$ is a smooth, bounded conformal factor. The smoothness and boundedness of $\omega$ are sufficient for the estimates in \autoref{rigidity conformally flat space}; an explicit expression for $\omega$ is not required.

Compared to Scheuer~\cite{scheuer_minkowski_2022}, the specialisation to $\mathbb{S}^n$ as base manifold is indeed a restriction to a more special setting. However, it is precisely this restriction that enables the stability results: local conformal flatness of the ambient space, which is needed for the rigidity argument in Step~2, holds for $\mathbb{S}^n$ but fails for general compact bases. In this sense, Theorems~\ref{theorem 1}-\ref{theorem 3} prove strictly more within a more special setting, and their conclusions have no analogue in the framework of \cite{scheuer_minkowski_2022}.
\end{rem}

\noindent We now impose the following structural conditions on the warping function, which we maintain throughout the paper.

\begin{assum}\label{assumption}
The warping function $\lambda \in C^{\infty}([a,b])$ satisfies:
 \begin{enumerate}
    \item $\lambda' > 0$;
    \item either $\lambda'' > 0$, or $\lambda'' \leq 0$ and $\partial_{r}\!\left(\tfrac{\lambda''}{\lambda}\right) \leq 0$.
    \end{enumerate}
\end{assum}

\noindent The condition $\lambda' > 0$ ensures that hypersurfaces can be represented as graphs over $\mathbb{S}^n$, and geometrically means that $M$ is expanding in the radial direction. The second condition controls the curvature: it guarantees that the Ricci curvature of $M$ satisfies the lower bounds needed for the monotonicity of the flow functionals (and hence for the Minkowski inequalities themselves), and it also enters the a priori estimates of \cite[Lemmas~3.3 and~3.5]{scheuer_minkowski_2022} that ensure uniform bounds on the principal curvatures and long-time existence of the flow. Both the RN-AdS and AdS-Schwarzschild manifolds satisfy \autoref{assumption}; we verify this explicitly in \autoref{RN AdS space}.

\medskip
We now state the five main stability theorems. The first three concern general warped products satisfying \autoref{assumption}; the last two concern the RN-AdS and AdS-Schwarzschild manifolds.

\subsection*{Stability in general warped product spaces}

Scheuer~\cite{scheuer_minkowski_2022} establishes three Minkowski inequalities for strictly convex graphs over $\mathbb{S}^n$ in warped product spaces. The first, under a Ricci curvature lower bound on $\mathbb{S}^n$ involving $\lambda$, gives
\begin{equation}
    \int_{\Sigma} H_1 + \frac{1}{n} \int_{\hat{\Sigma}} \overline{\mathrm{Ric}}(\partial_r, \partial_r) \geq \phi(|\Sigma|),
\end{equation}
where $\hat{\Sigma}$ is the region enclosed by $\Sigma$, $H_1$ is the normalised mean curvature, and $\phi$ depends on $|\Sigma|$. Equality holds if and only if $\Sigma$ is totally umbilic; if the relevant Ricci curvature quadratic forms are strictly positive, equality is attained only by radial slices. Specialising to $\mathbb{S}^n$ as base manifold, the Ricci curvature condition becomes the scalar inequality $\lambda'(r)^2 - \lambda(r)\lambda''(r) \leq 1$ (see \autoref{remark lower bounds on Ricci}).

\begin{thm}\label{theorem 1}
Let $n=2$. In addition to \autoref{assumption}, suppose that $M$ satisfies
\begin{equation}
    \lambda'(r)^2 - \lambda(r) \lambda''(r) \leq 1 \quad \textrm{for all } r.
\end{equation}
Let $\Sigma \subset M$ be a strictly convex graph over $\mathbb{S}^{n}$. Then there exists a constant $C= C(n, \|\omega\|_{\infty}, \|\nabla\omega\|_{\infty}, \mathrm{Vol}(\Sigma), \max_\Sigma|A|)$ such that
\begin{equation}\label{stability 1.6}
    \mathrm{dist}(\Sigma, S_{r}) \leq  C\left(\int_{\Sigma}H_{1} +\frac{1}{n}\int_{\hat{\Sigma}}\overline{\mathrm{Ric}}(\partial_{r}, \partial_{r})- \phi(|\Sigma|)\right)^{\!\frac{1}{2(n+1)}},
\end{equation}
where $S_r$ is a radial slice and $\omega$ is the conformal factor of the metric as in \autoref{remark conformal}.
\end{thm}

\begin{rem}
Constants $C$, $c$ may vary from line to line but remain independent of key geometric quantities unless specified otherwise.
\end{rem}

The second Minkowski inequality of \cite{scheuer_minkowski_2022} additionally requires convexity $\lambda'' \geq 0$ of the warping function and a Heintze-Karcher-type condition \eqref{Heintze-Karcher}, yielding
\begin{equation}
    \int_{\Sigma} H_1 + \frac{1}{n} \int_{\hat{\Sigma}} \overline{\mathrm{Ric}}(\partial_r, \partial_r) \geq \psi(|\hat{\Sigma}|),
\end{equation}
where $\psi$ depends on the enclosed volume $|\hat{\Sigma}|$, again with equality only for radial slices.

\begin{thm}\label{theorem 2}
Let $n=2$. In addition to \autoref{assumption}, suppose $\lambda'' \geq 0$ and
\begin{equation}
    \lambda'(r)^2 - \lambda(r) \lambda''(r) \leq 1 \quad \textrm{for all } r.
\end{equation}
Let $\Sigma \subset M$ be a strictly convex graph over $\mathbb{S}^{n}$ such that, for every such $\Sigma$,
\begin{equation}\label{Heintze-Karcher}
    \int_{\Sigma}\frac{\lambda'}{H_{1}} \geq \int_{\Sigma}u,
\end{equation}
where $u$ is the support function, with equality implying total umbilicity. Then there exists a constant $C= C(n, \|\omega\|_{\infty},\|\nabla\omega\|_{\infty},\mathrm{Vol}(\Sigma), \max_\Sigma|A|)$ such that
\begin{equation}\label{stability 1.8}
    \mathrm{dist}(\Sigma, S_{r}) \leq C\left(\int_{\Sigma}H_{1} +\frac{1}{n}\int_{\hat{\Sigma}}\overline{\mathrm{Ric}}(\partial_{r}, \partial_{r})- \psi(|\hat{\Sigma}|)\right)^{\!\frac{1}{2(n+1)}},
\end{equation}
where $S_r$ is a radial slice and $\omega$ is the conformal factor as in \autoref{remark conformal}.
\end{thm}

For the third inequality, the warping function takes the explicit form $\lambda(r) = \alpha\sinh r + \beta\cosh r$ with $\alpha \geq \beta \geq 0$, and the inequality reads \cite{scheuer_minkowski_2022}
\begin{equation}\label{minkowski ineq theorem 3}
    \int_{\Sigma} H_1 - |\hat{\Sigma}| \geq \phi(|\Sigma|),
\end{equation}
with equality only for radial slices.

\begin{thm}\label{theorem 3}
In addition to \autoref{assumption}, suppose
\begin{equation}
    \lambda(r)= \alpha\sinh{r}+\beta\cosh{r}, \quad \alpha\geq \beta\geq 0
\end{equation}
with at least one inequality strict, and $\alpha^{2}-\beta^{2} \leq 1$. Let $\Sigma \subset M$ be a strictly convex graph over $\mathbb{S}^{n}$. \\ Then there exists a constant $C= C(n, \|\omega\|_{\infty}, \|\nabla\omega\|_{\infty}, \mathrm{Vol}(\Sigma), \max_\Sigma|A|)$ such that
\begin{equation}\label{stability 1.10}
    \mathrm{dist}(\Sigma, S_{r})\leq C \left(\int_{\Sigma}H_{1} -|\hat{\Sigma}| - \phi(|\Sigma|)\right)^{\!\frac{1}{2(n+1)}},
\end{equation}
where $S_r$ is a radial slice and $\omega$ is the conformal factor as in \autoref{remark conformal}.
\end{thm}

\begin{rem}[On the restriction to $n = 2$ in Theorems \ref{theorem 1} and \ref{theorem 2}] This restriction is inherited from Scheuer~\cite{scheuer_minkowski_2022}. The flow~\eqref{flow_scheuer_1} has speed $F = H_2/H_1$, and long-time existence requires that $F$ vanish on the boundary of the positive cone $\Gamma_+$. This holds only when $n = 2$: for $n > 2$, it is possible for $H_2$ to remain positive even when one principal curvature tends to zero, which obstructs the convexity-preserving estimate. Since our stability analysis depends on the flow's long-time convergence, we inherit this dimensional restriction. \autoref{theorem 3} is free of this restriction because it uses the inverse mean curvature flow, for which convexity and convergence hold in all dimensions $n \geq 2$. We nevertheless keep $n$ as a parameter in the statements of Theorems~\ref{theorem 1} and~\ref{theorem 2} to make transparent that the monotonicity properties of the flow are valid in all dimensions; the restriction to $n=2$ enters only through the convergence theory.
\end{rem}

\subsection*{Stability in RN-AdS and AdS-Schwarzschild manifolds}\leavevmode

The Reissner-Nordstr\"om Anti-de Sitter manifold $M = [s_0, \infty) \times \mathbb{S}^{n-1}$ is equipped with the metric
\begin{equation}
    \bar{g} = \frac{1}{1 + \kappa^2 s^2 - 2m s^{2-n} + q^2 s^{4-2n}}\, ds^2 + s^2 g_{\mathbb{S}^{n-1}},
\end{equation}
where $m$ is the black hole mass, $q$ the charge, and $\kappa$ is related to the cosmological constant; the boundary $\{s_0\}\times\mathbb{S}^{n-1}$ is the horizon. The associated Minkowski inequality\cite{wang_minkowski-type_2015} for $\Sigma$ compact, mean-convex and star-shaped is
\begin{equation}
    \begin{split}
          \int_{\Sigma} fH\,d\mu - n(n-1)\kappa^{2}\int_{\Omega}f\,d\mathrm{vol} \geq&\; (n-1)f(\bar{s})^{2}\bar{s}^{n-2}|\mathbb{S}^{n-1}|\\
 -&(n-1)\kappa^{2}\bar{s}^{n}|\mathbb{S}^{n-1}| +(n-1)\kappa^{2}s_{0}^{n}|\mathbb{S}^{n-1}|,
    \end{split}
\end{equation}
where $\bar{s}=  \left(\frac{|\Sigma|}{|\mathbb{S}^{n-1}|}\right)^{1/(n-1)}$ is the areal radius, and \\$f(s) = \sqrt{1 + \kappa^2 s^2-2m s^{2-n} + q^2 s^{4-2n}}$. Equality holds if and only if $\Sigma$ is a coordinate sphere $\{s\}\times\mathbb{S}^{n-1}$. The AdS-Schwarzschild case ($q = 0$, $\kappa = 1$) reduces to \cite{brendle_minkowski_2016}
\begin{equation}
    \int_{\Sigma} f H \, d\mu - n(n-1) \int_{\Omega} f \, d\mathrm{vol} \geq (n-1) |\mathbb{S}^{n-1}|^{\frac{1}{n-1}} \left( |\Sigma|^{\frac{n-2}{n-1}} - |\partial M|^{\frac{n-2}{n-1}} \right).
\end{equation}

\begin{thm}\label{theorem 4}
Let $\Sigma$ be a compact mean-convex, star-shaped hypersurface in the RN-AdS manifold $(M, \bar{g})$, with $\Omega$ the region bounded by $\Sigma$ and the horizon, and let $f$ satisfy conditions (H1)-(H5) of \autoref{RN AdS space}. There exists a constant $C= C(n, \|\omega\|_{\infty}, \|\nabla\omega\|_{\infty},\mathrm{Vol}(\Sigma), \max_\Sigma|A|)$ such that
\begin{equation}
    \begin{split}
         \mathrm{dist}(\Sigma, S_{r}) \leq\; & C \bigg( \int_{\Sigma} fH\,d\mu - n(n-1)\kappa^{2}\int_{\Omega}f\,d\mathrm{vol} -(n-1)f(\bar{s})^{2}\bar{s}^{n-2}|\mathbb{S}^{n-1}|  \\
& +(n-1)\kappa^{2}\bar{s}^{n}|\mathbb{S}^{n-1}| -(n-1)\kappa^{2}s_{0}^{n}|\mathbb{S}^{n-1}|\bigg)^{\!\frac{1}{2n}},
    \end{split}
\end{equation}
where $\bar{s}=  \left(\frac{|\Sigma|}{|\mathbb{S}^{n-1}|}\right)^{1/(n-1)}$ and $S_r$ is a coordinate sphere.
\end{thm}

\begin{thm}\label{theorem 5}
Let $\Sigma$ be a compact mean-convex, star-shaped hypersurface in the AdS-Schwarzschild manifold $(M, \bar{g})$, with $\Omega$ the region bounded by $\Sigma$ and the horizon and $f=\sqrt{1-ms^{2-n}+s^{2}}$. There exists a constant $C= C(n, \|\omega\|_{\infty},\|\nabla\omega\|_{\infty},\mathrm{Vol}(\Sigma), \max_\Sigma|A|)$ such that
\begin{equation}\label{stability AdS}
    \begin{split}
           \mathrm{dist}(\Sigma, S_{r}) \leq\;  &C\bigg(\int_{\Sigma} fH\,d\mu -n(n-1)\int_{\Omega} f\,d\mathrm{vol}\\&- (n-1)|\mathbb{S}^{n-1}|^{\frac{1}{n-1}}\left(|\Sigma|^{\frac{n-2}{n-1}}- |\partial M|^{\frac{n-2}{n-1}}\right)\bigg)^{\!\frac{1}{2n}},
    \end{split}
\end{equation}
where $S_r$ is a coordinate sphere.
\end{thm}

\begin{rem}\label{remark relationship}
Theorems~\ref{theorem 4} and~\ref{theorem 5} hold in all dimensions $n \geq 2$ and are not special cases of Theorems~\ref{theorem 1}-\ref{theorem 3} for two reasons. First, the hypersurfaces here are not required to be convex, only mean-convex and star-shaped. Second, the warping function for RN-AdS and AdS-Schwarzschild is defined implicitly via the static potential $f(s)$, which involves the physical parameters $m$, $q$, $\kappa$ in a non-trivial way. Working directly with the defining equation
\begin{equation*}
    (\bar{\Delta} f)\bar{g} - \bar{D}^2 f + f\,\overline{\mathrm{Ric}} = (n-2)(n-1)q^2 f s^{4-2n} g_{\mathbb{S}^{n-1}}
\end{equation*}
and the structural conditions (H1)-(H5) yields stability estimates with explicit dependence on the physical parameters $(m,q,\kappa)$, providing bounds that are not accessible via the framework of Theorems~\ref{theorem 1}-\ref{theorem 3}. The dimensional restriction to $n=2$ in Theorems~\ref{theorem 1} and~\ref{theorem 2} also does not apply here; see the remark following \autoref{theorem 3}.
\end{rem}
\subsection*{Proof strategy}

Throughout this paper we use $(M, \bar{g})$ for the ambient manifold 
and $(\Sigma, g)$ for the hypersurface with induced metric $g$. Geometric quantities associated with $\bar{g}$ carry an overbar; those associated with $g$ are left unmarked. We write $h$ for the second fundamental form of $\Sigma$ and $A = (h^i_j) = (g^{ik}h_{kj})$ for the associated Weingarten operator.  Since indices are raised with $g$, the two have the same pointwise norm, and in particular $|\mathring{A}| = |\mathring{h}|$; we accordingly refer to either as the 
traceless second fundamental form. All five proofs follow the same two-step strategy.

\medskip
\noindent\textit{Step 1: Quantitative flow analysis.} Each inequality is proved in the literature by identifying suitable monotone quantities along a locally constrained inverse curvature flow. The flow exists for all time and converges smoothly to a radial slice; the monotonicity of these quantities encodes the inequality, and equality characterises totally umbilic hypersurfaces.  For stability, we make this convergence quantitative: by integrating the time-derivative of a suitable deficit functional and selecting a time at which the integrand is minimised, we extract a bound on the $L^p$-norm of the traceless second fundamental form $\mathring{A}$ in terms of the deficit $\epsilon$ in the inequality. Specifically, we obtain
\begin{equation*}
    \int_{\Sigma^\epsilon}|\mathring{A}|^{n+1} \leq C\sqrt{\epsilon},
\end{equation*}
where $\Sigma^\epsilon$ is a hypersurface along the flow that is already close to the initial surface $\Sigma$.

\medskip
\noindent\textit{Step 2: Rigidity in locally conformally flat spaces.} We prove a new rigidity theorem (\autoref{rigidity conformally flat space}) that, given a bound on $\|\mathring{A}\|_{L^p}$, concludes that the hypersurface is Hausdorff-close to a radial slice. This extends a result of De~Rosa-Gioffr\'e \cite{rosa_absence_2021} from the Euclidean setting to locally conformally flat ambient spaces. As noted in \autoref{remark conformal}, the warped product spaces $(a,b)\times\mathbb{S}^n$ are locally conformally flat, and it is this property that allows Step~2 to be carried out. Combining both steps via the triangle inequality completes the proof.

\medskip

The rigidity theorem that drives Step~2 of the proof strategy is the following. It is new and of independent interest beyond the stability results proved here.

\begin{thm}\label{rigidity conformally flat space}
Let $n \geq 2$, and let $\Sigma$ be a closed hypersurface embedded in a locally conformally flat manifold $(M^{n+1}, \bar{g})$ with conformal factor $\omega$. Suppose there exist constants $c_0 > 0$ and $p > n$ such that
\begin{equation}
    \|h\|_{L^p(\Sigma)} \leq c_0.
\end{equation}
Then there exists $\delta_{0}>0$ depending on $n,p,c_{0},\|\omega\|_\infty,\|\nabla\omega\|_\infty$ such that if
\begin{equation}
    \|\mathring{h}\|_{L^p(\Sigma)} \leq \delta_0,
\end{equation}
the Euclidean conformal image $\tilde{\Sigma}$ of $\Sigma$ can be written as a radial graph over a round sphere:
\begin{equation}
    \tilde{\Sigma} = \{ e^{f(x)} x \mid x \in \mathbb{S}^n \},
\end{equation}
where $f \in W^{2,p}(\mathbb{S}^n)$ satisfies
\begin{equation}
    \|f\|_{W^{2,p}(\mathbb{S}^n)} \leq C \|\mathring{h}\|_{L^p(\Sigma)},
\end{equation}
for a constant $C$ depending on $n, p, c_0, \|\omega\|_{\infty}, \|\nabla \omega\|_{\infty}$, and $\mathrm{Vol}(\Sigma)$. In particular, $\Sigma$ is Hausdorff-close to a totally umbilical hypersurface in $(M, \bar{g})$.
\end{thm}

\subsection*{Acknowledgements}
I thank Julian Scheuer, Federica Dragoni, Nicolas Dirr and John Harvey for their support and encouragement.

\section{Rigidity Result}\label{rigidity result}

In this section we prove \autoref{rigidity conformally flat space} by extending \cite[Theorem~1.3]{rosa_absence_2021} to locally conformally flat spaces.

\begin{definition}[Conformally Equivalent Manifolds]
Two Riemannian manifolds $(M^{n+1}, \bar{g})$ and $(M^{n+1}, \tilde{g})$ are \emph{conformally equivalent} if there exists a smooth function $\omega: M \to \mathbb{R}$ such that $\bar{g}_{ij} = e^{2\omega}\tilde{g}_{ij}$.
\end{definition}

\begin{definition}[Locally Conformally Flat Manifold]
A Riemannian manifold $(M^{n+1}, \bar{g})$ is \emph{locally conformally flat} if every point $p \in M$ has a neighbourhood $U$ on which $\bar{g}|_U$ is conformally equivalent to a flat metric.
\end{definition}

We record how the second fundamental form transforms under a conformal change $\bar{g}_{ij} = e^{2\omega}\tilde{g}_{ij}$ of the ambient metric. Denoting Euclidean quantities by a tilde:
\begin{equation}\label{conformal_changes_formulaes}
    \begin{split}
         g_{ij}&=e^{2\omega}\tilde{g}_{ij},\\
         \nu &= e^{-\omega} \tilde{\nu},\\
         h_{ij}e^{-\omega} &= \tilde{h}_{ij} + \omega_{\beta} \tilde{\nu}^{\beta} \tilde{g}_{ij},\\
         e^{\omega} h_{i}^{j} &= \tilde{h}_{i}^{j} + \omega_{\beta} \tilde{\nu}^{\beta} \delta_{i}^{j}, \\
         e^{\omega} H &= \tilde{H} + n \omega_{\beta} \tilde{\nu}^{\beta}.
    \end{split}
\end{equation}
Here $H$ is the mean curvature, $\nu$ is outward unit normal, and $h_{ij}$ is the second fundamental form.

Throughout the following proof, $|\cdot|_{g}$ and $\|\cdot\|_{L^p(\Sigma)}$ denote norms with respect to the induced metric $g$ on $\Sigma$, while $|\cdot|_{\tilde g}$ and $\|\cdot\|_{L^p(\tilde\Sigma)}$ denote norms with respect to the Euclidean metric $\tilde g$ on the conformal image $\tilde\Sigma$.

\begin{proof}[Proof of \autoref{rigidity conformally flat space}]
We apply \cite[Theorem~1.3]{rosa_absence_2021} in the Euclidean setting and interpret the conclusion back in $(M,\bar{g})$. To apply the Euclidean theorem we need
\begin{equation}
    \|\tilde{h}\|_{L^{p}(\tilde{\Sigma})} \leq c_1 \quad \text{and} \quad \mathrm{Vol}(\tilde{\Sigma}) = \mathrm{Vol}_{g_{\mathbb{S}^{n}}}(\mathbb{S}^{n}).
\end{equation}
The volume condition may be arranged by rescaling; see \autoref{volume constraint}.

Using \eqref{conformal_changes_formulaes},
\begin{equation}
    |\tilde{h}|_{\tilde g}^{2} =e^{2\omega}|h|_g^{2} -2(\omega_{\beta}\tilde{\nu}^{\beta})\tilde{H} -n(\omega_{\beta}\tilde{\nu}^{\beta})^{2}.
\end{equation}
Applying Cauchy--Schwarz with parameter $\epsilon > 0$,
\begin{equation}
    -2(\omega_\beta \tilde{\nu}^\beta)\tilde{H} \leq \epsilon (\omega_\beta \tilde{\nu}^\beta)^2 + \frac{\tilde{H}^2}{\epsilon},
\end{equation}
and choosing $\epsilon$ appropriately and using $\tilde{H}^2 \leq n|\tilde{h}|_{\tilde g}^{2}$, we find that $\|h\|_{L^p(\Sigma)} \leq c_0$ implies $\|\tilde{h}\|_{L^p(\tilde{\Sigma})} \leq c_1$ for some $c_1 = c_1(c_0, \|\omega\|_\infty, \|\nabla\omega\|_\infty)$.

For the traceless part, the relation $\mathring{h}_{ij}=h_{ij}-\frac{H}{n}g_{ij}$ gives
\begin{equation}
    |\mathring{h}|_g^{2}=e^{-2\omega}|\mathring{\tilde{h}}|_{\tilde g}^{2}.
\end{equation}
Under the conformal change the volume element transforms as $d\tilde{\mu} = e^{-n\omega}d\mu$, so
\begin{equation}
    \|\mathring{\tilde{h}}\|_{L^p(\tilde{\Sigma})}^p
= \int_{\Sigma} |\mathring{h}|_{g}^p \cdot e^{(p-n)\omega} \, d\mu,
\end{equation}
giving
\begin{equation}
    \|\mathring{\tilde{h}}\|_{L^p(\tilde{\Sigma})}
\leq e^{\frac{p-n}{p} \|\omega\|_\infty} \|\mathring{h}\|_{L^p(\Sigma)}.
\end{equation}
Hence if $\|\mathring{h}\|_{L^{p}(\Sigma)} \leq \delta_{0}$, then $\|\mathring{\tilde{h}}\|_{L^{p}(\tilde{\Sigma})} \leq \delta_{1}$ for $\delta_1 = \delta_1(\delta_0, \|\omega\|_\infty)$.

Theorem~1.3 of \cite{rosa_absence_2021} is purely Euclidean: it takes the bound $\|\tilde h\|_{L^p(\tilde\Sigma)}\leq c_1$ and returns a threshold $\delta_1$ and constant $C'$ depending only on $(n,p,c_1)$. Since $c_1 = c_1(c_0,\|\omega\|_\infty,\|\nabla\omega\|_\infty)$ as above, both $\delta_1$ and $C'$ inherit this dependence, and translating back through $\|\mathring{\tilde h}\|_{L^p(\tilde\Sigma)}\leq e^{\frac{p-n}{p}\|\omega\|_\infty}\|\mathring h\|_{L^p(\Sigma)}$ shows that $\delta_0$ and $C$ depend on $n,p,c_0,\|\omega\|_\infty,\|\nabla\omega\|_\infty$, as stated in the theorem.

Applying \cite[Theorem~1.3]{rosa_absence_2021}, there exists $c \in \mathbb{R}^{n+1}$ such that $\tilde{\Sigma} - c = \{ e^{f(x)} x \mid x \in \mathbb{S}^{n} \}$ with
\begin{equation}
    \|f\|_{W^{2,p}(\mathbb{S}^{n})} \leq C' \|\mathring{\tilde{h}}\|_{L^p(\tilde{\Sigma})} \leq C \|\mathring{h}\|_{L^p(\Sigma)}.
\end{equation}
This shows that the Euclidean conformal image $\tilde{\Sigma}$ of $\Sigma$ is close to a round sphere. Since the traceless second fundamental form is conformally invariant up to scaling, and totally umbilical hypersurfaces remain totally umbilical under conformal transformations, we conclude that the original hypersurface $\Sigma$ is Hausdorff-close to a totally umbilical hypersurface in $(M, \bar{g})$.

In model spaces such as the sphere or hyperbolic space, where totally umbilical hypersurfaces are geodesic spheres, this further implies that $\Sigma$ is close to a geodesic sphere.

\end{proof}

\begin{rem}\label{volume constraint}
We assume that the volume constraint $\mathrm{Vol}_{\tilde{g}}(\tilde{\Sigma}) = \mathrm{Vol}(\mathbb{S}^n)$ holds. We can achieve this by rescaling the ambient metric $\bar{g}$ by a constant factor $\lambda > 0$, chosen so that the Euclidean volume of the image $\tilde{\Sigma}$ becomes equal to that of the standard sphere. We then apply the Euclidean rigidity theorem in the rescaled setting, and scale back to the original metric afterward. This rescaling does not affect the validity of the assumptions on the second fundamental form or its traceless part, since under a conformal scaling the $L^p$ norms transform by a power of $\lambda$, and these transformations preserve the smallness conditions. The only change is that the constants in the conclusion will now depend on the initial volume of $\Sigma$ in the metric $\bar{g}$.

\end{rem}

\section{Preliminaries for Minkowski Inequalities in Warped Product Spaces}
\label{sec:preliminaries}

We outline the flow, notation, and monotonicity results from \cite{scheuer_minkowski_2022} that form the foundation for our proofs.

The ambient space is $M = (a, b) \times \mathbb{S}^{n}$ with metric $\bar{g} = dr^2 + \lambda^2(r) g_{\mathbb{S}^{n}}$. Evolving hypersurfaces $\Sigma_t$ are graphs over $\mathbb{S}^n$, parametrised by\\ $\Sigma_t = \{ (r(t,y), y) : y \in \mathbb{S}^n\}$. The flow equation is
\begin{equation}
    \frac{\partial}{\partial t}X(\xi, t)= \left( \frac{\lambda'(r)}{F(\kappa)} - u \right) \nu(\xi, t),
\end{equation}
where $\kappa = (\kappa_1,\dots,\kappa_n)$ are the principal curvatures, $F$ is a symmetric homogeneous degree-$1$ curvature function, and $u$ is the support function. Under suitable conditions on $\lambda$ and
\begin{equation}
    F=\frac{H_{2}}{H_{1}},
\end{equation}
the flow becomes
\begin{equation}\label{flow_scheuer_1}
    \frac{\partial}{\partial t}X = \left( \lambda'(r)\frac{H_{1}}{H_{2}} - u \right) \nu.
\end{equation}
It is shown in \cite{scheuer_minkowski_2022} that \eqref{flow_scheuer_1} exists for all time and converges smoothly to a radial slice, preserving convexity.

\begin{rem}\label{remark lower bounds on Ricci}
The results of \cite{scheuer_minkowski_2022} are formulated for warped products over a general compact base manifold $(S_0, \sigma)$; we specialise throughout to $S_0 = \mathbb{S}^n$ with $\sigma = g_{\mathbb{S}^n}$. A key assumption in \cite[Theorems~1.6 and~1.8]{scheuer_minkowski_2022} is
\begin{equation} \label{ricci-lower-bound}
    \widehat{\mathrm{Ric}} \geq (n-1)\left( \lambda'(r)^2 - \lambda''(r)\lambda(r) \right) \sigma \quad \text{for all } r,
\end{equation}
which ensures monotonicity of the flow functionals. For the round sphere $S_0 = \mathbb{S}^n$, where $\widehat{\mathrm{Ric}} = (n-1)g_{\mathbb{S}^n}$, this simplifies to the scalar condition
\begin{equation}
    \lambda'(r)^2 - \lambda(r) \lambda''(r) \leq 1 \quad \text{for all } r.
\end{equation}
For \cite[Theorem~1.10]{scheuer_minkowski_2022}, the condition $\widehat{\mathrm{Ric}} \geq (n-1)(\alpha^2-\beta^2)$ similarly becomes $\alpha^2-\beta^2 \leq 1$.
\end{rem}

We use the normalised elementary symmetric polynomials
\begin{equation}
    H_k = \frac{1}{\binom{n}{k}} \sum_{1 \leq i_1 < \dots < i_k \leq n} \kappa_{i_1} \dots \kappa_{i_k},
\end{equation}
and write $H = nH_1$ for the mean curvature. The Minkowski-type formula for warped product spaces is
\begin{equation}
    \int_{\Sigma} u H_1 = \int_{\Sigma} \lambda',
\end{equation}
and for $H_2$,
\begin{equation}
    \int_{\Sigma} u H_2 = \int_{\Sigma} \lambda'(r) H_1 - \frac{1}{n(n-1)} \int_{\Sigma} \overline{\mathrm{Ric}}(\nu, \nabla \Theta),
\end{equation}
where $\Theta$ is a primitive of $\lambda(r)$.

The key monotonicity properties along \eqref{flow_scheuer_1} are:
\begin{itemize}
\item $\partial_t |\Sigma_t| \geq 0$, with equality if and only if $\Sigma_t$ is totally umbilic. This ensures the flow increases surface area unless it is already a radial slice.
\item $\partial_t \!\left( \int_{\Sigma_t} H_1 + \tfrac{1}{n}\int_{\hat\Sigma_t} 
\overline{\mathrm{Ric}}(\partial_r,\partial_r) \right) \leq 0$, with equality if 
and only if $\Sigma_t$ is totally umbilic. This drives the hypersurface towards 
the equality configuration.
\end{itemize}
Together, these are the central tools used in \cite{scheuer_minkowski_2022} to 
prove the Minkowski inequalities; in the present work we make them quantitative 
to obtain control on the traceless second fundamental form in terms of the deficit.

\section{Stability of Minkowski-type Inequalities in Warped Product Spaces}
Before turning to the individual proofs, we record a monotonicity property of the comparison functions that is used throughout this section.

\begin{rem}\label{remark phi increasing}
The comparison functions $\phi$ and $\psi$ appearing in Theorems~\ref{theorem 1} and~\ref{theorem 2} are defined implicitly through the radial slices: $\phi$ via $W(S_r) = \phi(|S_r|)$ and $\psi$ via $W(S_r) = \psi(|\hat S_r|)$, where
\begin{equation}
    W(S_r) = \int_{S_r}H_1 + \tfrac{1}{n}\int_{\hat S_r}\overline{\mathrm{Ric}}(\partial_r,\partial_r).
\end{equation}
Since $H_1 = \lambda'/\lambda$ and $\overline{\mathrm{Ric}}(\partial_r,\partial_r) = -n\lambda''/\lambda$ on radial slices, a direct computation gives
\begin{equation}
    \frac{d}{dr}W(S_r) = (n-1)|\mathbb{S}^n|\lambda(r)^{n-2}\lambda'(r)^2 > 0,
\end{equation}
so that $W(S_r)$, together with $|S_r|$ and $|\hat S_r|$, is strictly increasing in $r$; the functions $\phi$ and $\psi$ are therefore strictly increasing, and $\phi^{-1}$, $\psi^{-1}$ are well-defined and increasing. This guarantees, in the proofs of Theorems~\ref{theorem 1} and~\ref{theorem 2} below, that $Q(t)$ is monotone and $\epsilon \geq 0$. The analogous functional $W_2(S_r)$ for \autoref{theorem 3} does not admit the same cancellation for general $\lambda$; its monotonicity is established separately in the proof of \autoref{theorem 3}, using the specific form $\lambda = \alpha\sinh r + \beta\cosh r$.
\end{rem}

\subsection{Proof of \autoref{theorem 1}}
\pf{
Define the deficit
\begin{equation}
\epsilon = \phi^{-1}\!\left(\int_{\Sigma}H_{1} +\frac{1}{n}\int_{\hat{\Sigma}}\overline{\mathrm{Ric}}(\partial_{r}, \partial_{r})\right)- |\Sigma| \geq 0,
\end{equation}
and set
\begin{equation}
    Q(t)= \phi^{-1}\!\left(\int_{\Sigma_{t}}H_{1} +\frac{1}{n}\int_{\hat{\Sigma}_t}\overline{\mathrm{Ric}}(\partial_{r}, \partial_{r}) \right)-|\Sigma_t|.
\end{equation}
Since
\begin{equation}
    \int_{\Sigma_t} H_1 + \frac{1}{n} \int_{\hat{\Sigma}_t} \overline{\textrm{Ric}}(\partial_{r}, \partial_{r})
\end{equation}
is monotonically decreasing $|\Sigma_t|$ is monotonically increasing along the flow \eqref{flow_scheuer_1}, where $|\Sigma_{t}|$ evolves as
\begin{equation}\label{area evolution}
    \frac{\partial}{\partial t}|\Sigma_t| = \int_{\Sigma_t} 
    \left(\frac{\lambda' H_1}{H_2} - u\right) H \, d\mu.
\end{equation} 
We take the derivative of $Q(t)$ and integrate over time:
\begin{equation}
    \int_{t=0}^{\infty} \frac{\partial}{\partial t} Q(t) = \int_{t=0}^{\infty} \frac{\partial}{\partial t}\left( \phi^{-1}\left(\int_{\Sigma_t}H_{1} + \frac{1}{n}\int_{\hat{\Sigma}_t}\overline{\textrm{Ric}}(\partial_{r}, \partial_{r})\right) - |\Sigma_t|\right).
\end{equation}
Integrating $\partial_t Q$ gives $Q(\infty) - Q(0) = -\epsilon$, and hence
\begin{equation}
\int_{0}^{\infty} \frac{\partial}{\partial t} |\Sigma_t|\, dt \leq \epsilon.
\end{equation}
Using 
\begin{equation}
H = nH_1,
\end{equation}
\begin{equation}
(n-1)H_2 = \frac{1}{n}H^2 - \frac{1}{n}|A|^2,
\end{equation}
and
\begin{equation}
|A|^2 = |\mathring{A}|^2 + \frac{1}{n}H^2,
\end{equation}
 we obtain
\begin{equation}
    n\left(\frac{\lambda' H_1^2}{H_2} - uH_1 \right) = n(\lambda' - uH_1) + \frac{\lambda'}{(n-1)H_2}|\mathring{A}|^2.
\end{equation}
Integrating over $\Sigma_t$ and applying the Minkowski-type formula
\begin{equation}
    \int_\Sigma uH_1 = \int_\Sigma \lambda',
\end{equation}
we arrive at
\begin{equation}
\int_{0}^{\infty}\int_{\Sigma_{t}} \frac{\lambda'|\mathring{A}|^2}{H_2} \leq (n-1) \epsilon.
\end{equation}
We get 
\begin{equation}
     \underset{t \in [0, \sqrt{\epsilon}]}{\textrm{min}}\int_{\Sigma_{t}} \frac{\lambda'|\mathring{A}|^2}{H_2} \leq (n-1) \sqrt{\epsilon}.
\end{equation}
Let $\Sigma^{\epsilon}$ be the hypersurface where the minimum is attained 
\begin{equation}
     \int_{\Sigma^{\epsilon}} \frac{\lambda'|\mathring{A}|^2}{H_2} \leq (n-1) \sqrt{\epsilon}.
\end{equation} Using the uniform bounds on the principal curvatures \cite[Lemma~3.5]{scheuer_minkowski_2022}, which depend on $\max_\Sigma|A|$, and the lower bound $\lambda' \geq 2\alpha > 0$ \cite[Lemma~3.3]{scheuer_minkowski_2022}, we obtain
\begin{equation}\label{tracelessnorm1}
    \begin{split}
         \int_{\Sigma^{\epsilon}}|\mathring{A}|^{n+1} & \leq C \int_{\Sigma^{\epsilon}}|\mathring{A}|^{2} \\
& \leq C \frac{\textrm{max}_{\Sigma^{\epsilon}}H_{2}}{\textrm{min}_{\Sigma^{\epsilon}}\lambda'}\int_{\Sigma^{\epsilon}} \frac{\lambda'|\mathring{A}|^2}{H_2} \\
& \leq C \sqrt{\epsilon}.
    \end{split}
\end{equation}
We also want to estimate the Hausdorff distance between $\Sigma_{t}$ and $\Sigma_{0}=\Sigma$. Let $X(\xi, 0)$ and $X(\xi,t)$ be two points in $\Sigma_{0}$ and $\Sigma_{t}$ respectively. Let $\gamma: [0,t] \to M$ be a curve defined as 
\begin{equation}
    \gamma(\tau)= X(\xi,\tau).
\end{equation}
Then we have  due to \cite[lemma 3.3]{scheuer_minkowski_2022} ,
\begin{equation}
    \begin{split}
        \mathrm{d}_{M}(X(\xi,0), X(\xi,t))
& \leq t\max_{[0,t]}\lvert \del_{\tau}\gamma \rvert \leq 2Ct.
    \end{split}
\end{equation}
From this, we get 
\begin{equation}\label{distance}
\mathrm{dist}(\Sigma_{t}, \Sigma) \leq Ct \leq C\sqrt{\epsilon} \quad \forall\, t \in [0,\sqrt\epsilon].
\end{equation}
Hence from \eqref{tracelessnorm1} and \eqref{distance}, we get 
\begin{equation}
    \int_{\Sigma^{\epsilon}} |\mathring{A}|^{n+1} \leq C\sqrt{\epsilon} \quad \textrm{and} \quad \textrm{dist}(\Sigma^{\epsilon}, \Sigma)\leq C\sqrt{\epsilon}.
\end{equation}
Applying \autoref{rigidity conformally flat space} to $\Sigma^\epsilon$ yields that $\Sigma^\epsilon$ is Hausdorff-close to a radial slice $S_r$, and the triangle inequality then gives the estimate \eqref{stability 1.6} for $\Sigma$.

}

\subsection{Proof of \autoref{theorem 2}}
\pf{
The proof follows the same structure as that of \autoref{theorem 1}, with the enclosed volume $|\hat\Sigma|$ playing the role of the surface area $|\Sigma|$. We define
\begin{equation}
\epsilon = \psi^{-1}\!\left(\int_{\Sigma}H_{1} +\frac{1}{n}\int_{\hat{\Sigma}}\overline{\mathrm{Ric}}(\partial_{r}, \partial_{r})\right)- |\hat{\Sigma}| \geq 0.
\end{equation}
Setting $Q(t) = \psi^{-1}\!\left(\int_{\Sigma_{t}}H_{1} +\frac{1}{n}\int_{\hat{\Sigma}_t}\overline{\mathrm{Ric}}(\partial_{r}, \partial_{r}) \right) - |\hat\Sigma_t|$ and integrating $\partial_t Q$ as before gives
\begin{equation}
\int_{0}^{\infty}\frac{\partial}{\partial t}|\hat{\Sigma}_{t}|\, dt \leq \epsilon.
\end{equation}
The evolution of $|\hat\Sigma_t|$ along \eqref{flow_scheuer_1} is
\begin{equation}
    \partial_t|\hat\Sigma_t| = \int_{\Sigma_t}(\tfrac{\lambda'H_1}{H_2} - u).
\end{equation}
Using the identity
\begin{equation}
    \frac{\lambda'H_{1}}{H_{2}} -u = \frac{H_{1}^{2}}{H_{2}} \left(\frac{\lambda'}{H_{1}}-u\right) + \frac{u}{n(n-1)H_{2}}|\mathring{A}|^{2},
\end{equation}
together with $H_1^2 \geq H_2$ from \cite[Lemma~3.5]{scheuer_minkowski_2022}, which depend on $\max_\Sigma|A|$ and the Heintze-Karcher condition \eqref{Heintze-Karcher}, we obtain
\begin{equation}
\frac{1}{n(n-1)}\int_{0}^{\infty} \int_{\Sigma_{t}} \frac{u}{H_{2}}|\mathring{A}|^{2} \leq \epsilon.
\end{equation}
Selecting $\Sigma^\epsilon = \Sigma_{t^\epsilon}$ for some $t^\epsilon \in [0, \sqrt\epsilon]$ as before and using the bounds on $H_2$, and $u \geq c > 0$ from \cite[Lemma~3.3, equation~2.8]{scheuer_minkowski_2022}, we get
\begin{equation}
    \int_{\Sigma^{\epsilon}}|\mathring{A}|^{n+1} \leq C \sqrt{\epsilon}.
\end{equation}
The estimate on $\mathrm{dist}(\Sigma^\epsilon, \Sigma)$ follows identically from \eqref{distance}, and the conclusion follows by applying \autoref{rigidity conformally flat space} and the triangle inequality.

}

\subsection{Proof of \autoref{theorem 3}}

\pf{
We denote by $S_{r}$ the radial slices in $M$ and introduce the functional 
$W_{2}(\Sigma)= \int_{\Sigma}H_{1}-|\hat{\Sigma}|$. Consider the inverse mean 
curvature flow
\begin{equation}\label{flow:theorem3}
    \frac{\partial}{\partial t}X(\xi, t)= \frac{1}{H}\nu(\xi, t).
\end{equation}
Along this flow we track the quantity
\begin{equation}\label{Wt1.10}
    W_{2}(\Sigma_{t})= \int_{\Sigma_{t}}H_{1}- |\hat{\Sigma}_{t}|,
\end{equation}
and to measure the deviation from radial symmetry we define
\begin{equation}
    Q(t)= W_{2}(\Sigma_{t})- \phi(|\Sigma_{t}|),
\end{equation}
where $\phi$ is determined implicitly by radial slices via the identity 
$W_{2}(S_{r})= \phi(|S_{r}|)$, so that $Q(t) \geq 0$ with equality precisely 
when $\Sigma_t$ is a radial slice.

As in \autoref{remark phi increasing}, we first verify that $\phi$ is well-defined here by checking that $W_2(S_r)$ and $|S_r|$ are strictly increasing in $r$. Since $S_r$ has induced metric 
$\lambda(r)^2 g_{\mathbb{S}^n}$ and constant principal curvatures 
$\lambda'(r)/\lambda(r)$,
\begin{equation}
    |S_r| = \lambda(r)^n |\mathbb{S}^n|, \qquad 
    W_2(S_r) = \lambda'(r)\lambda(r)^{n-1}|\mathbb{S}^n| - |\mathbb{S}^n|\int_a^r \lambda(\rho)^n \, d\rho.
\end{equation}
Hence $\frac{d}{dr}|S_r| = n\lambda(r)^{n-1}\lambda'(r)|\mathbb{S}^n| > 0$ by 
$\lambda' > 0$, and
\begin{equation}
    \frac{d}{dr}W_2(S_r) = |\mathbb{S}^n|\lambda(r)^{n-2}
    \left(\lambda''(r)\lambda(r) + (n-1)\lambda'(r)^2 - \lambda(r)^2\right).
\end{equation}
Since $\lambda = \alpha\sinh r + \beta \cosh r$ satisfies $\lambda'' = \lambda$, 
the bracket simplifies to $(n-1)\lambda'(r)^2$, so
\begin{equation}
    \frac{d}{dr}W_2(S_r) = (n-1)|\mathbb{S}^n|\lambda(r)^{n-2}\lambda'(r)^2 > 0.
\end{equation}
Thus $\phi$ is well-defined and strictly increasing.

We begin by computing the evolution of the mean curvature integral along 
\eqref{flow:theorem3}. Along the inverse mean curvature flow, $H$ evolves as
\begin{equation}\label{evolution H}
    \frac{\partial}{\partial t}H = -\Delta\!\left(\frac{1}{H}\right) 
    -\frac{1}{H}\left(|A|^{2}+ \overline{\mathrm{Ric}}(\nu, \nu)\right).
\end{equation}
Integrating this over $\Sigma_t$, and decomposing $|A|^2 = |\mathring{A}|^2 + \tfrac{1}{n}H^2$,
\begin{equation}\label{meancuravture1.10}
    \begin{split}
        \partial_{t}\int_{\Sigma_{t}}H
        &= \int_{\Sigma_{t}}H- \int_{\Sigma_{t}}\frac{1}{H} 
        \left(|A|^{2}+ \overline{\mathrm{Ric}}(\nu, \nu)\right)\\
        &= \frac{n-1}{n}\int_{\Sigma_{t}}H- \int_{\Sigma_{t}}\frac{|\mathring{A}|^{2}}{H}
        - \int_{\Sigma_{t}}\frac{\overline{\mathrm{Ric}}(\nu,\nu)}{H}.
    \end{split}
\end{equation}
To bound the last term, we recall from \cite[equation~2.9]{scheuer_minkowski_2022} 
that $\overline{\mathrm{Ric}}(\nu, \nu) \geq -n\lambda''/\lambda$. For the warping 
function $\lambda(r) = \alpha \sinh r + \beta \cosh r$, a direct computation gives 
$\lambda''/\lambda = 1$, and hence $\overline{\mathrm{Ric}}(\nu, \nu) \geq -n$. 
Substituting into \eqref{meancuravture1.10} yields
\begin{equation}
    \partial_{t}\int_{\Sigma_{t}}H \leq \frac{n-1}{n}\int_{\Sigma_{t}}H
    - \int_{\Sigma_{t}}\frac{|\mathring{A}|^{2}}{H} + n\int_{\Sigma_{t}}\frac{1}{H}.
\end{equation}
The enclosed volume evolves as $\partial_{t}|\hat{\Sigma}_{t}| = 
\int_{\Sigma_{t}}H^{-1}$. Differentiating the identity $W_2(S_r) = \phi(|S_r|)$ 
gives
\begin{equation}
    \phi'(|S_{r(t)}|)|S_{r(t)}| = \frac{n-1}{n} \left(W_{2}(S_{r(t)})
    + |\hat{S}_{r(t)}|\right).
\end{equation}
For a general mean-convex hypersurface $\Sigma_t$, we select at each time $t$ 
a corresponding radial slice $S_{r(t)}$ satisfying $|\Sigma_{t}|= |S_{r(t)}|$.

Setting $\Sigma = \Sigma_0$ and writing the initial deficit as
\begin{equation}
    \epsilon = \int_{\Sigma}H_{1}- |\hat{\Sigma}| - \phi(|\Sigma|) \geq 0,
\end{equation}
we differentiate $Q(t)$ explicitly and apply the evolution inequalities above:
\begin{equation}
    \begin{split}
        \frac{\partial}{\partial t}\left(\int_{\Sigma_{t}}H_{1}- |\hat{\Sigma}_{t}|
        - \phi(|\Sigma_{t}|) \right) \leq\; & \frac{n-1}{n^{2}}\int_{\Sigma_{t}}H
        - \frac{1}{n}\int_{\Sigma_{t}}\frac{|\mathring{A}|^{2}}{H}  \\
        &+\int_{\Sigma_{t}}\frac{1}{H}-\int_{\Sigma_{t}}\frac{1}{H}
        -\phi'(|\Sigma_{t}|)|\Sigma_{t}|.
    \end{split}
\end{equation}
Combining this with \eqref{Wt1.10} and using the identity above for 
$\phi'(|S_{r(t)}|)|S_{r(t)}|$,
\begin{equation}
    \begin{split}
        \frac{d}{dt}Q(t) \leq\; & \frac{n-1}{n}\left(W_{2}(\Sigma_{t})+ |\hat{\Sigma}_{t}|\right)
        - \frac{1}{n}\int_{\Sigma_{t}}\frac{1}{H}|\mathring{A}|^{2}-\phi'(|\Sigma_{t}|)|\Sigma_{t}|\\
        =\; & \frac{n-1}{n} \left(W_{2}(\Sigma_{t})-W_{2}(S_{r(t)})\right)
        +\frac{n-1}{n} \left(|\hat{\Sigma}_{t}|- |\hat{S}_{r(t)}|\right)
        \\&-\frac{1}{n}\int_{\Sigma_{t}}\frac{1}{H}|\mathring{A}|^{2}.
    \end{split}
\end{equation}
Since $|S_{r(t)}| = |\Sigma_t|$, the first term on the right is exactly 
$\frac{n-1}{n}Q(t)$. The second term, $|\hat{\Sigma}_t| - |\hat{S}_{r(t)}|$, is 
non-positive by the isoperimetric inequality \cite{guanVolumePreservingFlow2019}, since 
$\Sigma_t$ and $S_{r(t)}$ enclose equal areas and the sphere maximises enclosed 
volume among hypersurfaces of given area. Discarding this term, we obtain the 
differential inequality
\begin{equation}\label{gronwall Q}
    \frac{d}{dt}Q(t) \leq \frac{n-1}{n}Q(t) - \frac{1}{n}\int_{\Sigma_{t}}\frac{|\mathring{A}|^{2}}{H}.
\end{equation}
Multiplying \eqref{gronwall Q} by the integrating factor $e^{-\frac{n-1}{n}t}$ 
and integrating over $[0,T]$ for arbitrary $T>0$ gives
\begin{equation}
    \frac{1}{n}\int_{0}^{T} e^{-\frac{n-1}{n}s}\int_{\Sigma_{s}}\frac{|\mathring{A}|^{2}}{H}\,ds 
    \leq Q(0) - e^{-\frac{n-1}{n}T}Q(T) \leq \epsilon,
\end{equation}
where we used $Q(T) \geq 0$, which holds for every finite $T$ since $\Sigma_T$ 
is a strictly convex graph over $\mathbb{S}^n$ (convexity being preserved along the 
flow) and hence satisfies \eqref{minkowski ineq theorem 3} directly. Letting 
$T \to \infty$,
\begin{equation}
    \int_{0}^{\infty} e^{-\frac{n-1}{n}s}\int_{\Sigma_{s}}\frac{|\mathring{A}|^{2}}{H}\,ds \leq n\epsilon.
\end{equation}
Assuming without loss of generality that $\epsilon \leq 1$ and restricting to 
$t \in [0,\sqrt{\epsilon}]$, where $e^{-\frac{n-1}{n}t} \geq e^{-\frac{n-1}{n}}$, 
we conclude
\begin{equation}
    \underset{t \in [0,\sqrt{\epsilon}]}{\min} \int_{\Sigma_{t}} 
    \frac{|\mathring{A}|^{2}}{H}\,d\mu \leq C\sqrt{\epsilon}.
\end{equation}
Let $\Sigma^{\epsilon}$ be the hypersurface where the minimum on the left is attained
\begin{equation}
    \int_{\Sigma^{\epsilon}} \frac{|\mathring{A}|^{2}}{H}\,d\mu \leq C\sqrt{\epsilon}.
\end{equation}
Using the uniform bounds on the second fundamental form and mean curvature from 
\cite[Lemma~3.5]{scheuer_minkowski_2022}, which depend on $\max_\Sigma|A|$,
\begin{equation}
    \begin{split}
         \int_{\Sigma^{\epsilon}} |\mathring{A}|^{n+1} 
         &\leq C\int_{\Sigma^{\epsilon}}|\mathring{A}|^{2}\\
         & \leq C\, \max_{\Sigma^{\epsilon}}H \int_{\Sigma^{\epsilon}} 
         \frac{|\mathring{A}|^{2}}{H}\\
         & \leq C \sqrt{\epsilon}.
    \end{split}
\end{equation}
With this bound established, controlling $\mathrm{dist}(\Sigma^\epsilon, \Sigma)$ 
via the flow speed and applying \autoref{rigidity conformally flat space} via the 
triangle inequality proceeds exactly as in the proof of \autoref{theorem 1}.
}

\section{Stability of Minkowski-type Inequalities in RN-AdS and AdS-Schwarzschild Manifolds}\label{RN AdS space}

In this section we study the stability of Minkowski-type inequalities for hypersurfaces in the AdS-Schwarzschild and RN-AdS manifolds. We describe the geometry, recall the relevant Minkowski inequality, and prove Theorems~\ref{theorem 4} and~\ref{theorem 5} by working in a more general setting of which both are special cases.

\begin{rem}
In contrast to the previous sections, where the ambient manifolds are of dimension $n+1$, in this section all ambient spaces are $n$-dimensional, so hypersurfaces are of dimension $n-1$.
\end{rem}

\begin{defn}
The AdS-Schwarzschild manifold is $M= [s_{0}, \infty) \times \mathbb{S}^{n-1}$ with metric
\begin{equation}
    \bar{g}= \frac{1}{1-ms^{2-n}+s^{2}}\, ds\otimes ds +s^{2}g_{\mathbb{S}^{n-1}},
\end{equation}
where $m > 0$ and $s_0$ is the unique positive root of $1 + s_0^2 - ms_0^{2-n} = 0$. The boundary $\{s_0\}\times\mathbb{S}^{n-1}$ is the horizon.
\end{defn}

Setting $f = \sqrt{1-ms^{2-n}+s^{2}}$, the function $f$ is the static potential, satisfying
\begin{equation}
    (\bar{\Delta} f)g- \bar{D}^{2}f +f \overline{\mathrm{Ric}}=0,
\end{equation}
and $\bar\Delta f = nf$ (from the trace).

\begin{defn}
The Reissner-Nordstr\"om Anti-de Sitter manifold is $M = [s_0, \infty) \times \mathbb{S}^{n-1}$ with metric
\begin{equation}
    \bar{g} = \frac{1}{1 + \kappa^2 s^2 - 2m s^{2-n} + q^2 s^{4-2n}}\, ds^2 + s^2 g_{\mathbb{S}^{n-1}},
\end{equation}
where $m, q, \kappa > 0$ with $q < m$, $\kappa < \infty$. Setting $q = 0$ and $\kappa = 1$ recovers the AdS-Schwarzschild manifold. The static potential is
\begin{equation}
    f(s) = \sqrt{1 + \kappa^2 s^2 - 2m s^{2-n} + q^2 s^{4-2n}},
\end{equation}
satisfying
\begin{equation}
    (\bar{\Delta} f)\bar{g} - \bar{D}^2 f + f\overline{\mathrm{Ric}} = (n-2)(n-1) q^2 f s^{4-2n} g_{\mathbb{S}^{n-1}}.
\end{equation}
\end{defn}

By a change of variables \cite[Lemma~9]{wang_minkowski-type_2015}, both manifolds can be written in warped product form $\bar{g} = dr^2 + \lambda^2(r)g_{\mathbb{S}^{n-1}}$ with $\lambda(r) = s$, $ds/dr = f(s)$. In this form they are locally conformally flat, allowing us to apply \autoref{rigidity conformally flat space}.

\smallskip
We verify that both manifolds satisfy \autoref{assumption}. The condition $\lambda'(r)^2 - \lambda(r)\lambda''(r) \leq 1$ translates, in $(s,f)$ coordinates, to $f^2 - sff' \leq 1$. For AdS-Schwarzschild:
\begin{align*}
       f^2 - s f f' &= f^2 - \frac{s}{2}(f^2)'\\& = (1 + s^2 - ms^{2-n}) - (s^2 + \frac{n-2}{2}ms^{2-n})\\& = 1 - \frac{n}{2}ms^{2-n}\\
       & < 1
\end{align*}
since $m > 0$ and $s > 0$. For RN-AdS:
\begin{align*}
f^2 - sff' = 1 - nms^{2-n} + (n-1)q^2 s^{4-2n},
\end{align*}
which satisfies $f^2 - sff' \leq 1$ in the physical regime $q < m$ outside the horizon (and strictly as $s\to\infty$).

\smallskip
We prove stability in the following general setting, of which both RN-AdS and AdS-Schwarzschild are special cases. Consider $M=[s_{0}, \infty)\times \mathbb{S}^{n-1}$ with metric $\bar{g}= f^{-2}ds^{2}+s^{2}g_{\mathbb{S}^{n-1}}$, where $f$ is continuous on $[s_0,\infty)$ with $s_0 > 0$ and satisfies:
\begin{itemize}
    \item (H1)$f$ is differentiable and positive on $(s_0, \infty)$, and $f(s_0) = 0$.
    \item (H2) $1 + \kappa^2 s^2 - f^2 = 2ms^{2-n} + O(s^{4-2n})$.
    \item (H3) $f' > 0$.
    \item (H4) $f(f'^2 + f f'') + (n-3) \frac{f^2 f}{s} + (n-2)(1 - f^2)\frac{f}{s^{2}}\geq 0$.
    \item (H5) $P(x) = R(x^{1/(n-1)})x$, where $R(s) = 2f(s)f'(s) - 2\kappa^2 s$, is non-decreasing and concave.
\end{itemize}
These conditions are exactly those under which the argument of \cite{wang_minkowski-type_2015} goes through: (H1)-(H2) give a well-defined horizon and correct asymptotic Anti-de Sitter behaviour; (H3) ensures $C^1$ regularity for long-time existence; (H4) controls the monotonicity formula for $\int_{\Sigma_t}fH\,d\mu$; and (H5) ensures the key monotonicity of $Q(t)$ (defined in \eqref{Q def RNADS}) along the flow. When $P(x) = 2(n-2)m - 2(n-2)q^2 x^{-(n-2)/(n-1)}$, one recovers the RN-AdS manifold; setting $q = 0$ and $\kappa = 1$ gives AdS-Schwarzschild.

Under (H1)-(H5), the following Minkowski inequality holds \cite{wang_minkowski-type_2015}: for $\Sigma$ compact, mean-convex and star-shaped in $(M, \bar{g})$,
\begin{equation}
    \begin{split}
      \int_{\Sigma} fH\,d\mu- n(n-1)\kappa^{2}\int_{\Omega}f\,d\mathrm{vol} &\geq (n-1)f(\bar{s})^{2}\bar{s}^{n-2}|\mathbb{S}^{n-1}|\\
&-(n-1)\kappa^{2}\bar{s}^{n}|\mathbb{S}^{n-1}| +(n-1)\kappa^{2}s_{0}^{n}|\mathbb{S}^{n-1}|,
    \end{split}
\end{equation}
with equality if and only if $\Sigma = \{s\}\times\mathbb{S}^{n-1}$ for some $s \geq s_0$, where $s$ is the radial coordinate. The flow is the inverse mean curvature flow $\partial_t X = H^{-1}\nu$, which preserves the star-shaped and mean-convex properties, exists for all time \cite[Cor.~18]{wang_minkowski-type_2015}, and moves hypersurfaces outward, away from the horizon \cite[Lemma~12]{wang_minkowski-type_2015}. Rather than converging to a fixed hypersurface, $\Sigma_t$ becomes asymptotically umbilical as $t\to\infty$, which gives $\liminf_{t\to\infty}Q(t)\geq 0$ \cite[Prop.~22]{wang_minkowski-type_2015}; combined with the monotonicity of $Q(t)$ \cite[Lemma~23]{wang_minkowski-type_2015}, this yields $Q(0)\geq 0$.

\subsection{Proof of \autoref{theorem 4} and \autoref{theorem 5}}

We prove the general case; the stated theorems follow by substituting the 
specific functions $f$.

\begin{proof}
Define the deficit
\begin{equation}
    \begin{split}
           \epsilon = \int_{\Sigma} fH\,d\mu &- n(n-1)\kappa^{2}\int_{\Omega}f\,
           d\mathrm{vol} -(n-1)f(\bar{s})^{2}\bar{s}^{n-2}|\mathbb{S}^{n-1}| \\
&+(n-1)\kappa^{2}\bar{s}^{n}|\mathbb{S}^{n-1}| 
-(n-1)\kappa^{2}s_{0}^{n}|\mathbb{S}^{n-1}|.
    \end{split}
\end{equation}
To track the evolution of this quantity along the inverse mean curvature flow, 
we introduce
\begin{align*}
         W(t)= \int_{\Sigma_{t}} fH\,d\mu &- n(n-1)\kappa^{2}
            \int_{\Omega_{t}}f\,d\mathrm{vol} 
            -(n-1)f(\bar{s}_{t})^{2}|\Sigma_{t}|^{\frac{n-2}{n-1}}
            |\mathbb{S}^{n-1}|^{\frac{1}{n-1}} \\
&+(n-1)\kappa^{2}\bar{s}_{t}^{n}|\mathbb{S}^{n-1}|-(n-1)\kappa^{2}s_{0}^{n}|\mathbb{S}^{n-1}|,
\end{align*}
where $\bar{s}_{t}= \left(|\Sigma_{t}|/|\mathbb{S}^{n-1}|\right)^{1/(n-1)}$, 
and set 
\begin{equation} \label{Q def RNADS}
    Q(t)=|\Sigma_{t}|^{-(n-2)/(n-1)} W(t).
\end{equation} 
Since 
$\partial_t|\Sigma_t| = |\Sigma_t|$ along the inverse mean curvature flow, differentiating $Q(t)$ gives
\begin{align} 
    |\Sigma_{0}|^{\frac{n-2}{n-1}}\left(Q(\infty)-Q(0)\right) = &|\Sigma_{0}|^{\frac{n-2}{n-1}}\int_{t=0}^{\infty}\frac{dQ}{dt}\notag\\
=& |\Sigma_{0}|^{\frac{n-2}{n-1}} \int_{t=0}^{\infty}-\frac{n-2}{n-1}|\Sigma_{t}|^{-\frac{n-2}{n-1}} W(t)\notag\\
&+ |\Sigma_{0}|^{\frac{n-2}{n-1}} \int_{t=0}^{\infty} |\Sigma_{t}|^{-\frac{n-2}{n-1}} \frac{d}{dt}W(t). \label{Q(T) RNADS}
\end{align}
We now estimate $\tfrac{d}{dt}W(t)$ by treating each term separately:
\begin{align}\label{W(t)-1}
     \int_{0}^{\infty} |\Sigma_{t}|^{-\frac{n-2}{n-1}} \frac{d}{dt}W(t)\,dt
          = & \int_{0}^{\infty} |\Sigma_{t}|^{-\frac{n-2}{n-1}} \bigg( 
          \frac{d}{dt}\int_{\Sigma_{t}}fH\,d\mu \notag \\
 & -n(n-1)\kappa^{2} \frac{d}{dt} \int_{\Omega_{t}} f\,d\mathrm{vol}\notag\\ 
 & -(n-1)|\mathbb{S}^{n-1}|^{\frac{1}{n-1}} 
 \frac{d}{dt}\!\left(f(\bar{s}_{t})^{2}|\Sigma_{t}|^{\frac{n-2}{n-1}}\right)\notag\\
 & +(n-1)\kappa^{2} |\mathbb{S}^{n-1}|
 \frac{d}{dt}(\bar{s}_{t}^{n})\bigg)dt.
\end{align}
For the first term, using the evolution of $H$ along the inverse mean curvature 
flow, \eqref{evolution H}, and the identity $\Delta f = \bar\Delta f - \bar{D}^2f(\nu,\nu) - H\langle\bar\nabla f,\nu\rangle$, a direct computation gives
\begin{align*}
      \frac{d}{dt} \left( \int_{\Sigma_t} fH\, d\mu \right) 
        &= - \int_{\Sigma_t} f \Delta \left( \frac{1}{H} \right) d\mu 
        - \int_{\Sigma_t} \frac{f}{H} (|A|^2 + \overline{\mathrm{Ric}}(\nu,\nu))
        \, d\mu \\
&\quad + \int_{\Sigma_t} \left( \langle \bar{\nabla} f, \nu \rangle 
+ fH \right) d\mu \\
&= -\int_{\Sigma_{t}} \frac{1}{H}\Delta f\,d\mu
- \int_{\Sigma_t} \frac{f}{H} (|A|^2 + \overline{\mathrm{Ric}}(\nu,\nu))
\, d\mu \\
&\quad + \int_{\Sigma_t} \left( \langle \bar{\nabla} f, \nu \rangle 
+ fH \right) d\mu \\
&= - \int_{\Sigma_t} \frac{1}{H} \left( \bar{\Delta} f - \bar{D}^2 f(\nu,\nu) 
+ f\, \overline{\mathrm{Ric}}(\nu,\nu) \right) d\mu \\
&\quad - \int_{\Sigma_t} \frac{f}{H} |A|^2\, d\mu 
+ \int_{\Sigma_t} \left( 2 \langle \bar{\nabla} f, \nu \rangle 
+ fH \right) d\mu \\
&\leq -\int_{\Sigma_{t}} \frac{f}{H}|\mathring{A}|^{2} 
+ \int_{\Sigma_{t}} 2\langle\bar{\nabla}f, \nu\rangle 
+ \frac{n-2}{n-1}\int_{\Sigma_{t}} fH.
\end{align*}
For the volume term, we rely on the Heintze-Karcher type inequality \cite{brendle_constant_2013}
\begin{equation}
    (n-1)\int_{\Sigma_t} \frac{f}{H}\, d\mu \geq n \int_{\Omega_t} f\, d\mathrm{vol} + s_0^n |\mathbb{S}^{n-1}|.
\end{equation}
Applying this to the evolution of the enclosed volume along our flow gives
\begin{equation}
\frac{d}{dt} \int_{\Omega_t} f\, d\mathrm{vol} = \int_{\Sigma_t} 
\frac{f}{H}\, d\mu \geq \frac{n}{n-1} \int_{\Omega_t} f\, d\mathrm{vol} 
+ \frac{1}{n-1} s_0^n |\mathbb{S}^{n-1}|,
\end{equation}
and hence
\begin{equation}
    -n(n-1)\kappa^{2}\frac{d}{dt} \int_{\Omega_{t}} f\,d\mathrm{vol} 
    \leq -n^{2}\kappa^{2} \int_{\Omega_{t}} f\,d\mathrm{vol} 
    -n\kappa^{2}s_{0}^{n} |\mathbb{S}^{n-1}|.
\end{equation}
For the remaining terms, using $\tfrac{d}{dt}\bar{s}_{t}= 
\tfrac{1}{n-1}\bar{s}_{t}$,
\begin{equation}
    \begin{split}
         -(n-1)|\mathbb{S}^{n-1}| ^{\frac{1}{n-1}} \frac{d}{dt} \bigg(&f(\overline{s}_{t})^{2}|\Sigma_{t}|^{\frac{n-2}{n-1}}\bigg) \\=& -(n-1)|\mathbb{S}^{n-1}|\frac{d}{dt} \left(f(\overline{s}_{t})^{2}\overline{s}_{t}^{n-2}\right)\\
=& -(n-1)|\mathbb{S}^{n-1}| \bigg(\frac{2}{n-1}f(\overline{s}_{t})f'(\overline{s}_{t})\overline{s}_{t}^{n-1} \\&+ \frac{n-2}{n-1}\left(f(\overline{s}_{t})\right)^{2}\overline{s}_{t}^{n-2}\bigg)\\
=& -|\mathbb{S}^{n-1}| \bigg(2f(\overline{s}_{t})f'(\overline{s}_{t})\overline{s}_{t}^{n-1}+ (n-2)\left(f(\overline{s}_{t})\right)^{2}\overline{s}_{t}^{n-2}\bigg).
    \end{split}
\end{equation}
and
\begin{equation}
    (n-1)\kappa^{2} |\mathbb{S}^{n-1}| \frac{d}{dt}(\bar{s}_{t}^{n}) 
    = n\kappa^{2} |\mathbb{S}^{n-1}|\overline{s}_{t}^{n}. 
\end{equation}
Substituting all these estimates into \eqref{W(t)-1} yields
\begin{equation}\label{W(t)-2}
    \begin{split}
            \int_{t=0}^{\infty} |\Sigma_{t}|^{\frac{n-2}{n-1}}\frac{d}{dt}W(t) \leq &\int_{t=0}^{\infty} |\Sigma_{t}|^{-\frac{n-2}{n-1}}\bigg( -\int_{\Sigma_{t}}\frac{f}{H} |\mathring{A}|^{2}+ \int_{\Sigma_{t}}2\langle\overline{\nabla}f, \nu\rangle \\ &+ \frac{n-2}{n-1}\int_{\Sigma_{t}}fH-n^{2}\kappa^{2}\int_{\Omega_{t}} fd\textrm{vol}\\
& -n\kappa^{2}s_{0}^{n}|\mathbb{S}^{n-1}|\\&- |\mathbb{S}^{n-1}| \left(2f(\overline{s}_{t})f'(\overline{s}_{t})\overline{s}_{t}^{n-1}+ (n-2)\left(f(\overline{s}_{t})\right)^{2}\overline{s}_{t}^{n-2}\right)\\
&+n\kappa^{2}|\mathbb{S}^{n-1}|\overline{s}_{t}^{n} \bigg).
    \end{split}
\end{equation}
To handle the gradient term, we apply the identities from 
\cite[Lemma~23]{wang_minkowski-type_2015}:
\begin{equation}
       \int_{\Sigma_{t}} \langle\bar{\nabla}f, \nu\rangle\, d\mu
       = \int_{\mathbb{S}^{n-1}}ff's^{n-1}\big|_{s_{t}(w)}
       \,d\mathrm{vol}_{\mathbb{S}^{n-1}},
\end{equation}
\begin{equation}
    \begin{split}
            2\int_{\mathbb{S}^{n-1}}ff's^{n-1}\big|_{s_{t}(w)}
            \,d\mathrm{vol}_{\mathbb{S}^{n-1}}
            &= 2f(\bar{s}_{t})f'(\bar{s}_{t})s_{t}^{n-1}|\mathbb{S}^{n-1}|
            -2\kappa^{2}\bar{s}_{t}^{n}|\mathbb{S}^{n-1}|\\
            &\quad+2\kappa^{2}\int_{\mathbb{S}^{n-1}}s_{t}^{n}(w)
            \,d\mathrm{vol}_{\mathbb{S}^{n-1}}, 
    \end{split}
\end{equation}
and
\begin{equation}
    \int_{\mathbb{S}^{n-1}}{s_{t}^{n}(w)}\,d\mathrm{vol}_{\mathbb{S}^{n-1}} 
    = n\int_{\Omega_{t}} f\,d\mathrm{vol}+ |\mathbb{S}^{n-1}|s_{0}^{n}.
\end{equation}
Combining these three identities,
\begin{equation}
    \begin{split}
          2\int_{\Sigma_{t}} \langle\overline{\nabla}f, \nu\rangle d\mu &= 2f(\bar{s}_{t})f'(\bar{s}_{t})s_{t}^{n-1}|\mathbb{S}^{n-1}|-2\kappa^{2}\bar{s}_{t}^{n}|\mathbb{S}^{n-1}|\\&+2\kappa^{2}\int_{\mathbb{S}^{n-1}}s_{t}^{n}(w)d\textrm{vol}_{\mathbb{S}^{n-1}} + 2\kappa^{2}n\int_{\Omega_{t}} fd\textrm{vol}\\&+ 2\kappa^{2}s_{0}^{n}|\mathbb{S}^{n-1}|.
    \end{split}
\end{equation}

Substituting this into \eqref{W(t)-2}, all terms except the traceless part 
cancel, giving
\begin{align}
       \int_{t=0}^{\infty} |\Sigma_{t}|^{-\frac{n-2}{n-1}} \frac{d}{dt}W(t) &\leq \int_{t=0}^{\infty}|\Sigma_{t}|^{-\frac{n-2}{n-1}} \bigg(\int_{\Sigma_{t}}\frac{f}{H} |\mathring{A}|^{2} +(n-2)\kappa^{2}|\mathbb{
 S}^{n-1}|\bar{s}_{t}^{n}\\&-n(n-2)\kappa^{2}\int_{\Omega_{t}}fd\textrm{vol}-(n-2)\kappa^{2}s_{0}^{n}|\mathbb{S}^{n-1}|\\&+\frac{n-2}{n-1}\int_{\Sigma_{t}}fHd\mu-(n-2)\left(f(\overline{s}_{t})\right)^{2}\overline{s}_{t}^{n-2}\bigg).
\end{align}

Inserting this back into the integrated evolution of $Q(t)$ from $0$ to $\infty$, the remaining terms cancel precisely, leaving:
\begin{equation}
       |\Sigma_{0}|^{\frac{n-2}{n-1}}\left(Q(\infty) - Q(0)\right) \leq - |\Sigma_{0}|^{\frac{n-2}{n-1}}\int_{0}^{\infty}|\Sigma_{t}|^{-\frac{n-2}{n-1}} \int_{\Sigma_{t}}\frac{f}{H}|\mathring{A}|^{2}\,d\mu\,dt.
\end{equation}
By the definition of our initial deficit, we have $|\Sigma_0|^{\frac{n-2}{n-1}}Q(0) = \epsilon$. Furthermore, because the flow surfaces become asymptotically umbilical as $t \to \infty$, we know that $\liminf_{t\to\infty}Q(t)\geq 0$. Therefore, the left-hand side is bounded below by $-\epsilon$:
\begin{equation}
       -\epsilon \leq |\Sigma_{0}|^{\frac{n-2}{n-1}}\left(Q(\infty) - Q(0)\right).
\end{equation}
Combining these inequalities and multiplying by $-1$, we obtain the desired bound on the traceless second fundamental form:
\begin{equation}
       |\Sigma_{0}|^{\frac{n-2}{n-1}}\int_{0}^{\infty}|\Sigma_{t}|^{-\frac{n-2}{n-1}} \int_{\Sigma_{t}}\frac{f}{H}|\mathring{A}|^{2}\,d\mu\,dt \leq \epsilon.
\end{equation}
Taking the minimum over $t \in [0, \sqrt{\epsilon}]$,
\begin{equation}
    \underset{t \in [0,\sqrt{\epsilon}]}{\min}\; |\Sigma_{0}|^{\frac{n-2}{n-1}} 
    |\Sigma_{t}|^{-\frac{n-2}{n-1}} \int_{\Sigma_{t}} \frac{f}{H} |\mathring{A}|^{2}
    \,d\mu \leq \sqrt{\epsilon}. 
\end{equation}
Let $\Sigma_{l}$ be the hypersurface where this minimum is attained, so that
\begin{equation}
    \left(\frac{|\Sigma_{0}|}{|\Sigma_{l}|} \right)^{\!\frac{n-2}{n-1}} 
    \int_{\Sigma_{l}} \frac{f}{H} |\mathring{A}|^{2}\,d\mu \leq \sqrt{\epsilon}.
\end{equation}
Since we are interested in the stability regime, where the deficit is small, we assume without loss of generality that $\epsilon \leq 1$. Then $l \in [0,\sqrt{\epsilon}] \subseteq [0,1]$, and since $|\Sigma_{l}| = e^{l} |\Sigma_{0}|$,
\begin{equation}
    \frac{|\Sigma_{l}|}{|\Sigma_{0}|} = e^{l}\leq e^{\sqrt{\epsilon}}\leq e.
\end{equation}
Hence the ratio $|\Sigma_0|/|\Sigma_l|$ is bounded uniformly, and using the 
uniform bounds on $H$ from \cite[Cor.~16 and Prop.~17]{wang_minkowski-type_2015}, which depend on $\max_\Sigma|A|$,
\begin{equation}
    \begin{split}
         \int_{\Sigma_{l}} |\mathring{A}|^{n} &\leq C \int_{\Sigma_{l}} 
         |\mathring{A}|^{2}\\
& \leq C \frac{\max_{\Sigma_{l}} H}{\min_{\Sigma_{l}} f} 
\int_{\Sigma_{l}} \frac{f}{H} |\mathring{A}|^{2}\,d\mu.
    \end{split}
\end{equation}
It remains to bound $\min f$ from below. Since $f=0$ only at the horizon $s=s_0$ and $f$ is strictly increasing in $s$ by (H1) and (H3), it suffices to bound $\min_{\Sigma_t} s$ away from $s_0$. The inverse mean curvature flow moves hypersurfaces outward, satisfying $\min_{\Sigma_t} s \geq e^{t/(n-1)}\min_{\Sigma_0} s$ \cite[Lemma~12]{wang_minkowski-type_2015}. On the time interval $t\in[0,\sqrt\epsilon]$ this is bounded below by the fixed positive constant $\min_{\Sigma_0} s$, giving a uniform constant $c > 0$ such that $f \geq c$ throughout the flow. Therefore,
\begin{equation}
    \int_{\Sigma_{l}} |\mathring{A}|^{n} \leq C \sqrt{\epsilon}.
\end{equation}
With this bound on $\|\mathring{A}\|$ established, controlling 
$\mathrm{dist}(\Sigma_l, \Sigma)$ via the flow speed and applying 
\autoref{rigidity conformally flat space} via the triangle inequality proceeds 
exactly as in the proof of \autoref{theorem 1}.

\end{proof}
\section*{Statements and Declarations}

\noindent \textbf{Funding.} Prachi Sahjwani was partially supported by a DTP Programme of EPSRC (Project Reference EP/T517951/1; sub-project \enquote{Stability in physical systems governed by curvature quantities}, Project Reference 2601534) and by a UKRI Future Leaders Fellowship [grant number MR/W01176X/1; PI: J.~Harvey].

\noindent \textbf{Competing interests.} The author has no competing interests to declare that are relevant to the content of this article.

\noindent \textbf{Data availability.} On behalf of all authors, the 
corresponding author states that the manuscript has no associated data.

\end{document}